\documentclass[11pt]{article}
\usepackage{amsmath, amsthm, amssymb, graphicx, txfonts, bm}
\usepackage{accents}
\usepackage{soul}
\usepackage{times}
\usepackage{tikz} 
\usepackage{pgf}
\usepackage{graphics}
\usepackage{mathrsfs}
\usepackage{framed}
\usepackage{comment}

\newtheorem{theorem}{Theorem}[section]
\newtheorem{lemma}[theorem]{Lemma}

\newtheorem{corollary}[theorem]{Corollary}

\numberwithin{equation}{section}

\usetikzlibrary{positioning}

\DeclareMathOperator{\divergenz}{div}

\newcommand{\dt}{\,dt}

\newcommand{\Rone}{\mathbb{R}}
\newcommand{\Rn}{\mathbb{R}^n}

\newcommand{\Comone}{\mathbb{C}}

\newcommand{\CR}{\mathscr{C}_{R,d}^n}
\newcommand{\Rnone}{\mathbb{R}^{n+1}}

\begin{document}

\title{Motion by Volume Preserving Mean Curvature Flow Near Cylinders}

\author{David Hartley\footnote{Affiliation: Monash University, \texttt{david.hartley@monash.edu}}}

\date{}

\maketitle

\begin{abstract}
Center manifold analysis can be used in order to investigate the stability of the stationary solutions of various PDEs. This can be done by considering the PDE as an ODE between certain Banach spaces and linearising about the stationary solution. Here we investigate the volume preserving mean curvature flow using such a technique. We will consider surfaces with boundary contained within two parallel planes such that the surface meets these planes orthogonally. The stationary solution we consider is a cylinder. We will find that for initial surfaces that are sufficiently close to a cylinder the flow will exist for all time and converge to a cylinder exponentially. In particular, we show that there exists global solutions to the flow that converge to a cylinder, which are initially non-axially symmetric. A similar case where the initial surfaces are compact without boundary has previously been investigated by Escher and Simonett \cite{Escher98A}.

\end{abstract}

\section{Introduction}
Given a sufficiently smooth, compact initial hypersurface $\Omega_0=\bm{X}_0(M^n)\subset\Rnone$ we are interested in finding a family of embeddings $\bm{X}:M^n\times[0,T)\rightarrow\Rnone$ such that
\begin{equation}
\frac{\partial\bm{X}}{\partial t}=\left(h-H\right)\bm{\nu}_{\Omega_t},\ \ \ \bm{X}(\cdot,0)=\bm{X}_0,\ \ \ h=\frac{1}{\left|\Omega_t\right|}\int_{M^n} H\,d\mu_t,
\label{FullFlow}
\end{equation}
where $H=\sum_{i=1}^n\kappa_i$ is the mean curvature and $\kappa_i$ are the principal curvatures of the surface $\Omega_t=\bm{X}(M^n,t)=\bm{X}_t(M^n)$, $\bm{\nu}_{\Omega_t}$ and $\,d\mu_t$ are a choice of unit normal and the induced measure of $\Omega_t$ respectively. Note that the flow (\ref{FullFlow}) decreases surface area while preserving volume (see Section \ref{NotPre}).

For convex hypersurfaces without boundary it is known that a unique solution exists for all time and that the hypersurfaces converge to a sphere as $t\rightarrow\infty$ \cite{Huisken87}. Other results include average mean convex surfaces with initially small traceless second fundamental form converging to spheres \cite{Li09} and surfaces that are graphs over spheres with a height function close to zero, in a certain function space, converging to spheres \cite{Escher98A}. The case where the initial surface has a boundary has also been studied. In this case it is assumed that $\Omega_0$ is smoothly embedded in the domain
\begin{equation*}
W=\left\{\bm{x}\in\Rnone:0<x_{n+1}<d\right\},
\end{equation*}
with $d>0$ and $\partial \Omega_0\subset\partial W$. The volume enclosed by $\Omega_0$ and $W$ will be labelled $V$ and $\bm{\nu}_{\Omega_t}$ is chosen to point outward. The boundary conditions for the flow are then that $\Omega_0$ intersects $\partial W$ orthogonally. Assuming $\Omega_0$ to be rotationally symmetric it was proved in \cite{Athanassenas97} that the flow exists for all time and converges to a cylinder in $W$ of volume $V$ under the assumption
\begin{equation}
\left|\Omega_0\right|\leq\frac{V}{d}.
\label{MariaAssump}
\end{equation}
This constraint is needed to ensure that the solution never touches the axis of rotation, so that no singularities develop. Short time existence of the flow can be proven using standard parabolic techniques for initially smooth surfaces, however in this paper we are concerned with a larger class of initial surfaces so need a more general existence theorem (see Section \ref{NotPre}).

The main result of this paper is
\begin{theorem}
If $\Omega_0$ is a graph over the cylinder $\CR$, such that $R>\frac{d\sqrt{n-1}}{\pi}$ and with height function sufficiently small in $h^{1+\beta}_{\frac{\partial}{\partial z}}\left(\CR\right)$, $0<\beta<1$, then its flow by (\ref{FullFlow}) exists for all time and converges exponentially fast to a cylinder as $t\rightarrow\infty$, with respect to the $C^k$ topology for any $k\in\mathbb{N}$.
\label{MainR}
\end{theorem}
An important point to note is that the limiting cylinder may not be the same as the cylinder $\CR$ thus the cylinders are stable critical points of the PDE (\ref{FullFlow}), but not necessarily asymptotically stable. An immediate consequence of this theorem is that there are non-axially symmetric surfaces that converge to a cylinder under the flow (\ref{FullFlow}), since there are non-axially symmetric surfaces in any $h^{1+\alpha}_{\frac{\partial}{\partial z}}$ neighbourhood of a cylinder. Lastly note that the condition $R>\frac{d\sqrt{n-1}}{\pi}$ is of a similar nature to the condition (\ref{MariaAssump}); in fact for a cylinder (\ref{MariaAssump}) reduces to $R\geq nd$. In view of this, when considering the flow of axially-symmetric surfaces that are close to a cylinder we reproduce the results of \cite{Athanassenas97} under a weaker condition. The condition is also sharp, in the sense that if $R=\frac{d\sqrt{n-1}}{\pi}$ then there is a one-parameter family of unduloids, satisfying the boundary conditions, such that $\CR$ is the only cylinder in the family, i.e., there are stationary solutions to the flow arbitarily close to $\CR$ that are not cylinders.

In Section \ref{NotPre} of this paper we convert the flow (\ref{FullFlow}) to a PDE for the graph function and also introduce the function spaces and notation that will be used throughout the paper. The section ends with some important theorems about the flow. In Section \ref{Stab} we consider the problem as an ODE on Banach spaces and determine the spectrum of the linearisation of the speed. This then allows us to find a center manifold of the flow, which is seen to consist locally entirely of cylinders. Lastly exponential attractivity to this manifold is proven.

The author would like to thank Maria Athanassenas for her initial suggestion of working with cylinders, support, advice and help in preparing this paper, and also Todd Oliynyk for his advice and encouragement in dealing with center manifolds. Finally to Monash University and the School of Mathematical Sciences for their support.

\section{Notation and Preliminaries}
\label{NotPre}
In this paper we consider $M^n=\CR$, a given cylinder of radius $R$ and length $d$, and $\bm{X}_0(\bm{p})=\bm{p}+\rho_0(\bm{p})\bm{\nu}_{\CR}(\bm{p})$, $\bm{p}\in\CR$, so that $\bm{X}_0$ is a graph over $\CR$. The volume form on such a graph surface will be denoted by $\,d\mu_{\rho}$ and we let $\mu_{\rho}$ be the function such that $\,d\mu_{\rho}=\mu_{\rho}\,d\mu_0$. We use the notation $\bm{q}\in\mathscr{S}_R^{n-1}$, the $(n-1)$-dimensional sphere of radius $R$, and $z\in[0,d]$, so that $(\bm{q},z)\in\CR$. Proceeding as Escher and Simonett \cite{Escher98A} and converting the flow to an evolution equation for the height function $\rho:\CR\times[0,T)\rightarrow\Rone$ it is known that up to a tangential diffeomorphism the flow (\ref{FullFlow}) is equivalent to 
\begin{equation}
\frac{\partial\rho}{\partial t}=\sqrt{1+\tilde{g}^{ij}_{\rho}\nabla_{i}\rho\nabla_j\rho}\left(h_{\rho}-H_{\rho}\right),\ \ \ \ \left.\frac{\partial\rho}{\partial z}\right|_{z=0,d}=0,\ \ \ \ \rho(\cdot,0)=\rho_0,
\label{GraphFlow}
\end{equation}
where $H_{\rho}=H_{\rho}(\cdot,t)$ is the mean curvature of the surface defined by $\rho(\cdot,t)$  and $\tilde{g}^{ij}_{\rho}$ denotes the inverse metric on $\mathscr{C}_{R+\rho,d}^n$ (see \cite{Bode07}).

The graph functions $\rho$ are chosen in the little-H\"older spaces, which are defined for $\alpha\in(0,1)$, $k\in\mathbb{N}$, an open set $U\subset \Rn$ and a multi-index $\beta=(\beta_1,\ldots,\beta_n)$ with $|\beta|=\sum_{i=1}^n\beta_i$ as follows (see \cite{Lunardi95}):
\begin{align*}
&h^{\alpha}\left(\bar{U}\right)=\left\{\rho\in C^{\alpha}\left(\bar{U}\right):\lim_{r\rightarrow0}\sup_{x,y\in \bar{U},\ 0<|x-y|<r}\frac{|\rho(x)-\rho(y)|}{|x-y|^{\alpha}}=0\right\},\\
&h^{k+\alpha}\left(\bar{U}\right)=\left\{\rho\in C^{k+\alpha}\left(\bar{U}\right): D^{\beta}\rho\in h^{\alpha}\left(\bar{U}\right)\textrm{ for all } \beta,\ |\beta|=k\right\},
\end{align*}
where $D$ is the derivative operator on $\Rn$ and $C^{\alpha}$, $C^{k+\alpha}$ are the H\"older spaces:
\begin{align*}
&C^{\alpha}\left(\bar{U}\right)=\left\{\rho\in C\left(\bar{U}\right):\sup_{x,y\in \bar{U},x\neq y}\frac{|\rho(x)-\rho(y)|}{|x-y|^{\alpha}}<\infty\right\},\\
&C^{k+\alpha}\left(\bar{U}\right)=\left\{\rho\in C^{k}\left(\bar{U}\right): D^{\beta}\rho\in C^{\alpha}\left(\bar{U}\right)\textrm{ for all }  \beta,\ |\beta|=k\right\}.
\end{align*}
The norms on the little-H\"older spaces are the same as those on the H\"older spaces,
\begin{equation*}
\|\rho\|_{h^{k+\alpha}\left(\bar{U}\right)}=\|\rho\|_{C^{k+\alpha}\left(\bar{U}\right)}=\|\rho\|_{C^k\left(\bar{U}\right)}+\sum_{|\beta|=k}\sup_{x,y\in \bar{U},x\neq y}\frac{|D^{\beta}\rho(x)-D^{\beta}\rho(y)|}{|x-y|^{\alpha}},
\end{equation*}
where
\begin{equation*}
\|\rho\|_{C^k\left(\bar{U}\right)}=\sum_{|\beta|\leq k}\sup_{x\in\bar{U}}\left|D^{\beta}\rho(x)\right|.
\end{equation*}
The little-H\"older spaces can be extended to a manifold by means of the atlas of $\CR$. In addition, it is known that the little-H\"older spaces are the continuous interpolation spaces between the $C^k$ spaces (see \cite{Lunardi95}),
\begin{equation*}
h^{k\theta}\left(\CR\right)=\left(C\left(\CR\right),C^k\left(\CR\right)\right)_{\theta},
\end{equation*}
where $(\cdot,\cdot)$ is an interpolation functor for each $\theta\in(0,1)$ and is defined for $Y\subset X$ as follows:
\begin{align*}
(X,Y)_{\theta}=\left\{x\in X:\lim_{t\to0}t^{-\theta}K(t,x,X,Y)=0\right\},\ K(t,x,X,Y)=\inf_{a\in Y}\left(\left\|x-a\right\|_{X}+t\left\|a\right\|_Y\right).
\end{align*}
The interpolation space $h^{\alpha}\left(\CR\right)$ can also be viewed as the closure of the smooth functions on the cylinder, $C^{\infty}\left(\CR\right)$, inside the space of H\"older continuous functions $C^{\alpha}\left(\CR\right)$ (see \cite{Lunardi95}). To simplify the notation involving the boundary conditions, let $Y\left(\CR\right)$ represent any function space on the cylinder and set
\begin{equation*}
Y_{\frac{\partial}{\partial z}}\left(\CR\right)=\left\{\rho\in Y\left(\CR\right):\left.\frac{\partial\rho}{\partial z}\right|_{z=0,d}=0\right\}.
\end{equation*}
With this notation we also have that
\begin{equation}
\left(h^{\alpha}\left(\CR\right),h_{\frac{\partial}{\partial z}}^{2+\alpha}\left(\CR\right)\right)_{\theta}=\left\{\begin{array}{ll} h^{\alpha+2\theta}\left(\CR\right) & \theta<\frac{1-\alpha}{2}\\ h_{\frac{\partial}{\partial z}}^{\alpha+2\theta}\left(\CR\right) & \theta>\frac{1-\alpha}{2},\theta\neq1-\frac{\alpha}{2}\end{array}\right.,
\label{Interp}
\end{equation}
where $\alpha\in(0,1)$ - a result that follows from the Reiteration Theorem (\cite{Lunardi95} Theorem 1.2.15) along with the characterisation of interpolation spaces in Chapter 3 of \cite{Lunardi95}.

\renewcommand{\labelenumi}{(\roman{enumi})}
For an operator between function spaces $G:Y\rightarrow \tilde{Y}$ we denote the Fr\'echet derivative by $\partial G$. A linear operator, $A:Y\subset X\rightarrow X$, is called \textit{sectorial} if there exist $\theta\in\left(\frac{\pi}{2},\pi\right)$ and $\omega\in\Rone$ and $M>0$ such that
\begin{enumerate}
\item $\rho(A)\supset S_{\theta,\omega}=\{\lambda\in\Comone:\lambda\neq\omega,||\arg(\lambda-\omega)|<\theta\}$,
\item $\|R(\lambda,A)\|_{L\left(X\right)}\leq\frac{M}{|\lambda-\omega|}\textrm{ for all }\lambda\in S_{\theta,\omega}$,
\end{enumerate}
here $\rho(A)$  is the resolvent set, $R(\lambda,A)=(\lambda I-A)^{-1}$ is the resolvent operator and $\|\cdot\|_{L\left(X\right)}$ is the standard linear operator norm (see \cite{Lunardi95}).

The following are some facts about the flow.
\begin{lemma}
For an initially smooth, compact hypersurface with boundary contained in $\partial W$, the flow (\ref{FullFlow}) together with Neumann boundary condition:
\begin{enumerate}
\item  preserves the enclosed volume $V$,
\item  decreases the surface area $\left|\Omega_t\right|$.
\end{enumerate}
\label{VolAre}
\end{lemma}

\proof
Both parts are covered in \cite{Athanassenas97}. Using the divergence theorem and setting $E_t\subset W$ to be the $(n+1)$-dimensional space with boundary $\Omega_t$ we have the first part of the lemma as follows
\begin{equation*}
\frac{dV}{\dt}=\int_{E_t}\divergenz \frac{\partial\bm{X}}{\partial t}\,d\bm{x}=\int_{\partial E_t}\frac{\partial\bm{X}}{\partial t}\cdot\bm{\nu}_{\partial E_t}\,dS=\int_{M^n} \left(h-H\right)\,d\mu_t=0,
\end{equation*}
where the integral over the part of the boundary contained in $\partial W$ vanishes due to the boundary conditions. The second statement is through a direct calculation of the evolution equation for the area
\begin{equation*}
\frac{d\left|\Omega_t\right|}{\dt}=-\int_{M^n}(H-h)^2\,d\mu_t.
\end{equation*}
$\Box$

Note that Lemma \ref{VolAre} is also true for hypersurfaces without boundary. We now give a general short-time existence theorem

\begin{theorem}
For $0<\beta<1$ and a given hypersurface $\Omega_0$ of class $h^{1+\beta}$, the flow (\ref{FullFlow}) has a unique classical solution $\left\{\Omega_t:t\in[0,T)\right\}$ for some $T>0$. Each hypersurface $\Omega_t$ is of class $C^{\infty}$ for $t\in(0,T)$. Moreover, the mapping $t\rightarrow\Omega_t$ is continuous on $[0,T)$ with respect to the $h^{1+\beta}$-topology and smooth on $(0,T)$ with respect to the $C^{\infty}$-topology.
\end{theorem}

\proof
This result follows from the main result by Escher and Simonett in \cite{Escher98A}. It is clear that the arguments in their proof still hold with $h^{2+\alpha}$, $h^{1+\beta}$ and $h^{1+\beta_0}$ replaced by $h_{\frac{\partial}{\partial z}}^{2+\alpha}$, $h_{\frac{\partial}{\partial z}}^{1+\beta}$ and $h_{\frac{\partial}{\partial z}}^{1+\beta_0}$ respectively due to equation (\ref{Interp}).
$\Box$

\section{Stability around Cylinders}
\label{Stab}

The flow in equation (\ref{GraphFlow}) can be considered as an ordinary differential equation between Banach spaces. Set $0<\alpha<\beta<1$ and define
\begin{equation*}
G:h^{2+\alpha}_{\frac{\partial}{\partial z}}\left(\CR\right)\rightarrow h^{\alpha}\left(\CR\right),\ \ \ G(\rho):=L_{\rho}\left(h_{\rho}-H_{\rho}\right),\ \ \ L_{\rho}:= \sqrt{1+\tilde{g}^{ij}\nabla_i\rho\nabla_j\rho}.
\end{equation*}
The flow (\ref{GraphFlow}) is then rewritten as
\begin{equation}
\rho'(t)=G(\rho(t)),\ \ \ \ 0<t<T,\ \ \ \ \ \ \rho(0)=\rho_0\in h^{1+\beta}_{\frac{\partial}{\partial z}}\left(\CR\right).
\label{GraphODE}
\end{equation}

\begin{lemma} For the linearisation, $\partial G(0)$, of $G$ it holds
\begin{equation*}
\partial G(0)u = \left(\frac{n-1}{R^2}+\Delta_{\CR}\right)u-\frac{n-1}{R^2}\fint_{\CR}u\,d\mu_0,
\end{equation*}
for $u\in h^{2+\alpha}_{\frac{\partial}{\partial z}}\left(\CR\right)$.
\label{LinOp2}
\end{lemma}

\proof
Firstly note that $L_0=1$ and that $\left.\partial L_{\rho}\right|_{\rho=0}=0$. Hence for $u\in h^{2+\alpha}_{\frac{\partial}{\partial z}}\left(\CR\right)$
\begin{align*}
\left.\partial h_{\rho}\right|_{\rho=0}u &= \left.\partial\left(\frac{1}{\int_{\CR}\mu_{\rho}\,d\mu_0}\int_{\CR}H_{\rho}\mu_{\rho}\,d\mu_0\right)\right|_{\rho=0}u\\
&=\frac{1}{\left(\int_{\CR}\,d\mu_0\right)^2}\left(\int_{\CR}\,d\mu_0\left.\partial\left(\int_{\CR}H_{\rho}\mu_{\rho}\,d\mu_0\right)\right|_{\rho=0}u\right.\\
&\hspace{4cm}\left. -\int_{\CR}H_{0}\,d\mu_0\left.\partial\left(\int_{\CR}\mu_{\rho}\,d\mu_0\right)\right|_{\rho=0}u\right)\\
&=\frac{1}{\left|\CR\right\|}\left(\int_{\CR}\left.\partial H_{\rho}\right|_{\rho=0}u+H_{0}\left.\partial\mu_{\rho}\right|_{\rho=0}u\,d\mu_0 -H_{0}\int_{\CR}\left.\partial\mu_{\rho}\right|_{\rho=0}u\,d\mu_0\right)\\
&=\fint_{\CR}\left.\partial H_{\rho}\right|_{\rho=0}u\,d\mu_0.
\end{align*}
It was shown in \cite{Escher98B} that
\begin{equation*}
\left.\partial H_{\rho}\right|_{\rho=0}=-\left(\frac{n-1}{R^2}+\Delta_{\CR}\right),
\end{equation*}
so combining these results gives, for $u\in h^{2+\alpha}_{\frac{\partial}{\partial z}}\left(\CR\right)$,
\begin{equation}
\partial G(0)u = \left(\frac{n-1}{R^2}+\Delta_{\CR}\right)u-\fint_{\CR}\left(\frac{n-1}{R^2}+\Delta_{\CR}\right)u\,d\mu_0.
\label{LinOp1}
\end{equation}
The divergence theorem and the boundary conditions give the result.
$\Box$

As we are considering the flow locally about $\rho=0$ it is convenient to rewrite (\ref{GraphODE}) highlighting the dominant linear part
\begin{equation}
\rho'(t)=\partial G(0)\rho(t)+\tilde{G}\left(\rho(t)\right),\ \ \ \ \tilde{G}\left(u\right):=G\left(u\right)-\partial G(0)u.
\label{GraphLin}
\end{equation}
\begin{lemma}
If $R>\frac{d\sqrt{n-1}}{\pi}$ then the spectrum $\sigma\left(\partial G(0)\right)$ of $\partial G(0)$ consists of a sequence of isolated non-positive eigenvalues where the multiplicity of the $0$ eigenvalue is $n+1$.
\label{Spectrum}
\end{lemma}

\proof
First note that $\partial G(0):h^{2+\alpha}_{\frac{\partial}{\partial z}}\left(\CR\right)\rightarrow h^{\alpha}\left(\CR\right)$ and since $h^{2+\alpha}_{\frac{\partial}{\partial z}}\left(\CR\right)$ is compactly embedded in $h^{\alpha}\left(\CR\right)$ the spectrum consists solely of eigenvalues. To characterise the spectrum of $\partial G(0)$ we first look at the spectrum of
\begin{equation*}
\tilde{A}:=\frac{n-1}{R^2}+\Delta_{\CR}.
\end{equation*}
In the following we deal with the natural complexification without distinguishing it notationally. The operator $\tilde{A}$ is self adjoint with respect to the $L_2$ inner product on $h_{\frac{\partial}{\partial z}}^{2+\alpha}\left(\CR\right)$, since the boundary terms disappear for functions in $h^{2+\alpha}_{\frac{\partial}{\partial z}}\left(\CR\right)$ when integrating by parts. This means that the spectrum consists of real numbers. For an eigenfunction $v$ of $\tilde{A}$ assume that $v(\bm{q},z)=X(\bm{q})Z(z)$, so using $\Delta_{\CR}=\Delta_{\mathscr{S}_R^{n-1}}+\frac{\partial^2}{\partial z^2}$:
\begin{equation*}
\left(\frac{n-1}{R^2}+\Delta_{\mathscr{S}_R^{n-1}}+\frac{\partial^2}{\partial z^2}\right)X(\bm{q})Z(z)=\lambda X(\bm{q})Z(z).
\end{equation*}
So, after expanding,
\begin{equation*}
Z(z)\Delta_{\mathscr{S}_R^{n-1}}X(\bm{q})+X(\bm{q})Z''(z)+\left(\frac{n-1}{R^2}-\lambda\right) X(\bm{q})Z(z)=0.
\end{equation*}
Therefore
\begin{equation*}
-\frac{\Delta_{\mathscr{S}_R^{n-1}}X(\bm{q})}{X(\bm{q})}=\frac{Z''(z)}{Z(z)}+\left(\frac{n-1}{R^2}-\lambda\right).
\end{equation*}
Setting both sides equal to a constant $\xi\in\Rone$
\begin{equation*}
\Delta_{\mathscr{S}_R^{n-1}}X(\bm{q})=-\xi X(\bm{q}).
\end{equation*}
The eigenvalues of the Laplacian on $\mathscr{S}_R^{n-1}$ are well known to be $\frac{-l(l+n-2)}{R^2}$ for $l\in\mathbb{N}\cup\{0\}$, with eigenfunctions the spherical harmonics of order $l$, $X_{l,p}(\bm{q})=Y_{l,p}(\bm{q})$, $1\leq p\leq M_l$, where
\begin{equation*}
M_l={}^{l+n-1}C_{n-1}-{}^{l+n-3}C_{n-1}.
\end{equation*}
Therefore $\xi_l=\frac{l(l+n-2)}{R^2}$ so the other separation equation becomes
\begin{equation*}
Z''(z)=\left(\lambda+\frac{l(l+n-2)-(n-1)}{R^2}\right)Z(z),\ \ \ Z'(0)=Z'(d)=0,
\end{equation*}
where we included the boundary conditions to ensure that $v\in h^{2+\alpha}_{\frac{\partial}{\partial z}}\left(\CR\right)$. The eigenvalues and eigenfunctions to the second separation equation are given by
\begin{equation*}
\lambda_{l,m}=-\left(\frac{m^2\pi^2}{d^2}+\frac{l(l+n-2)-(n-1)}{R^2}\right),\ \ \ m,l\in\mathbb{N}\cup\{0\}
\end{equation*}
and
\begin{equation*}
Z_m(z)=\cos\left(\frac{m\pi z}{d}\right).
\end{equation*}
So $\lambda_{l,m}$ are the eigenvalues of $\tilde{A}$ with multiplicity $M_l$ and eigenfunctions spanned by
\begin{equation*}
v_{l,p,m}(\bm{q},z)=\cos\left(\frac{m\pi z}{d}\right)Y_{l,p}(\bm{q}),\ \ \ 1\leq p\leq M_l.
\end{equation*}
Now since the spherical harmonics are dense in the continuous functions on $\mathscr{S}_R^{n-1}$ and the trigonometric functions are dense in the continuous functions on $[0,d]$, we have found every eigenfunction of $\tilde{A}$.

Returning to the spectrum of $\partial G(0)$, $v_{0,1,0}=1$ is still an eigenfunction but with eigenvalue $\lambda_{0,0}=0$. The operator $\partial G(0)$ is also self adjoint with respect to the $L_2$ inner product on $h_{\frac{\partial}{\partial z}}^{2+\alpha}\left(\CR\right)$. To see this, consider $u,w\in h^{2+\alpha}_{\frac{\partial}{\partial z}}\left(\CR\right)$, then
\begin{align*}
\int_{\CR}\left(\partial G(0)u\right)w\,d\mu_0 & =\int_{\CR}\left(\tilde{A}u-\frac{n-1}{R^2}\fint_{\CR}u\,d\mu_0\right)w\,d\mu_0\\
&=\int_{\CR}\left(\tilde{A}u\right)w\,d\mu_0-\frac{n-1}{R^2}\fint_{\CR}u\,d\mu_0\int_{\CR}w\,d\mu_0\\
&=\int_{\CR}u\left(\tilde{A}w\right)\,d\mu_0-\frac{n-1}{R^2}\int_{\CR}u\,d\mu_0\fint_{\CR}w\,d\mu_0\\
&=\int_{\CR}u\left(\tilde{A}w-\frac{n-1}{R^2}\fint_{\CR}w\,d\mu_0\right)\,d\mu_0\\
&=\int_{\CR}u\left(\partial G(0)w\right)\,d\mu_0.
\end{align*}
Therefore we need only consider other eigenfunctions, orthogonal to $v_{0,1,0}=1$, in order to characterise the spectrum. This means that for an eigenfunction $u$
\begin{equation*}
\int_{\CR}u\,d\mu_0=0,
\end{equation*}
hence by Lemma \ref{LinOp2}, $\partial G(0)u=\tilde{A}u$. The remaining eigenfunctions of $\partial G(0)$ are then the remaining eigenfunctions of $\tilde{A}$, with the same eigenvalues. Therefore for $m,l\in\mathbb{N}\cup\{0\}$ and $m^2+l^2\geq1$ the eigenvalues are
\begin{equation*}
\lambda_{l,m}=-\left(\frac{m^2\pi^2}{d^2}+\frac{l(l+n-2)-(n-1)}{R^2}\right),
\end{equation*}
with eigenfunctions
\begin{equation*}
v_{l,p,m}(\bm{q},z)=\cos\left(\frac{m\pi z}{d}\right)Y_{l,p}(\bm{q}),\ \ \ 1\leq p\leq M_l.
\end{equation*}
It is clear from these expressions that $\lambda_{1,0}=0$, and we relabel the constant eigenfunction as $v_{1,0,0}=1$ for convenience, so zero is an eigenvalue of multiplicity $M_1+1=n+1$. If $m\in\mathbb{N}$ we have $\lambda_{1,m}<0$. Also if $l\geq2$ and $m\in\mathbb{N}\cup{0}$ then $\lambda_{l,m}<0$. The last case is for $l=0$ and $m\geq1$ where the eigenvalues are bounded above by $\lambda_{0,1}$:
\begin{equation*}
\lambda_{0,m}=-\left(\frac{m^2\pi^2}{d^2}-\frac{n-1}{R^2}\right)\leq-\left(\frac{\pi^2}{d^2}-\frac{n-1}{R^2}\right)=\lambda_{0,1}.
\end{equation*}
These eigenvalues are all negative under the condition $R>\frac{d\sqrt{n-1}}{\pi}$ given in the statement of the theorem.
$\Box$

In what follows, we set $P$ to be the projection from $h^{\alpha}\left(\CR\right)$ onto the $\lambda=0$ eigenspace given by
\begin{equation*}
Pu:=\sum_{p=0}^{n}\left\langle u,v_{1,p,0}\right\rangle v_{1,p,0},
\end{equation*}
where we use $\left\langle \cdot,\cdot\right\rangle$ to denote the $L_2$ inner product on $h_{\frac{\partial}{\partial z}}^{2+\alpha}\left(\CR\right)$. Because $\partial G(0)$ is self adjoint with respect to this inner product, it is clear that $P\partial G(0)u=\partial G(0)Pu=0$ for every $u\in h^{2+\alpha}_{\frac{\partial}{\partial z}}\left(\CR\right)$. Due to this $h^{2+\alpha}_{\frac{\partial}{\partial z}}\left(\CR\right)$ can be split into the subspaces $X^c=P\left(h^{2+\alpha}_{\frac{\partial}{\partial z}}\left(\CR\right)\right)$ and $X^s=(I-P)\left(h^{2+\alpha}_{\frac{\partial}{\partial z}}\left(\CR\right)\right)$, called the \textit{center subspace} and \textit{stable subspace} respectively. We are now in a position to apply Theorem 9.2.2 in \cite{Lunardi95}.

\begin{theorem}
For any $k\in\mathbb{N}$, there is a function $\gamma\in C^{k-1}\left(X^c,X^s\right)$, such that $\gamma^{(k-1)}$ is Lipschitz continuous, $\gamma(0)=\partial\gamma(0)=0$ and $\mathcal{M}^c=\textrm{graph}(\gamma)$ is a locally invariant manifold of equation (\ref{GraphLin}) with dimension $n+1$ provided $R>\frac{d\sqrt{n-1}}{\pi}$.
\end{theorem}

Note that by \textit{locally invariant} it is meant that there exists a neighbourhood of zero in $\Lambda\subset X^c$ such that if $\rho_0\in\textrm{graph}\left(\gamma|_{\Lambda}\right)$ then the solution to (\ref{GraphLin}) is in $\textrm{graph}\left(\gamma|_{\Lambda}\right)$ for all time or until $P\rho(t)\notin\Lambda$.

\proof
In order to know that the assumptions of Theorem 9.2.2 in \cite{Lunardi95} are satisfied in our case, we use Equation (\ref{Interp}) to show that for any $\alpha_0$ satisfying $0<\alpha_0<\alpha<1$ the following holds
\begin{equation*}
\left(h^{\alpha_0}\left(\CR\right),h^{2+\alpha_0}_{\frac{\partial}{\partial z}}\left(\CR\right)\right)_{\frac{\alpha-\alpha_0}{2}}=h^{\alpha}\left(\CR\right).
\end{equation*}
The assumptions are therefore satisfied by proving there exists a neighbourhood of 0, $O\subset h_{\frac{\partial}{\partial z}}^{2+\alpha}\left(\CR\right)$,  such that $\tilde{G}\in C^{\infty}\left(O,h^{\alpha}\left(\CR\right)\right)$ and, for all $\rho\in O$, $\partial G(\rho)$ is the part in $h^{\alpha}\left(\CR\right)$ of a sectorial operator $A_{\rho}:h_{\frac{\partial}{\partial z}}^{2+\alpha_0}\left(\CR\right)\rightarrow h^{\alpha_0}\left(\CR\right)$.
It was shown in \cite{Bode07}, \cite{Escher98B} that the mean curvature can be written in the form 
\begin{equation*}
H_{\rho}=\frac{1}{L_{\rho}}\left(\frac{n-1}{R+\rho}+Q(\rho)\rho\right),
\end{equation*}
where $Q(\rho)$ is a second order uniformly elliptic operator that smoothly depends on $\rho$ and its first derivative. This means that $\tilde{G}$ is a smooth operator inside a neighbourhood, $O_1$, where if $\rho\in O_1$ we have  $\rho(\bm{p})>-R$ for all $\bm{p}\in\CR$.

As $Q(0)=-\Delta_{\CR}$ we have $\left.\partial H_{\rho}\right|_{\rho=0}=-\left(\frac{n-1}{R^2}+\Delta_{\CR}\right)$, this operator is uniformly elliptic, hence its negative is sectorial. Now the operator
\begin{align*}
&A_0:h_{\frac{\partial}{\partial z}}^{2+\alpha_0}\left(\CR\right)\rightarrow h^{\alpha_0}\left(\CR\right),\\
&A_0u=\left(\frac{n-1}{R^2}+\Delta_{\CR}\right)u-\frac{n-1}{R^2}\fint_{\CR}u\,d\mu_0,
\end{align*}
can be seen to be sectorial: using the definition of sectorial there exists $M>0$, $\theta\in\left(\frac{\pi}{2},\pi\right)$ and $\omega\in\Rone$ such that
\begin{equation*}
\left|\lambda-\omega\right|\left\|u\right\|_{h^{2+\alpha_0}_{\frac{\partial}{\partial z}}\left(\CR\right)}\leq M\left\|\left(\lambda+\left.\partial H_{\rho}\right|_{\rho=0}\right)u\right\|_{h^{\alpha_0}\left(\CR\right)},
\end{equation*}
for all $\lambda\in S_{\theta,\omega}$ and $u\in h^{2+\alpha_0}_{\frac{\partial}{\partial z}}\left(\CR\right)$. Therefore
\begin{align*}
\left\|(\lambda-A_0)u\right\|_{h^{\alpha_0}\left(\CR\right)}&=\left\|\left(\lambda-\left(\frac{n-1}{R^2}+\Delta_{\CR}\right)\right)u-\left(-\frac{n-1}{R^2}\fint_{\CR}u\,d\mu_0\right)\right\|_{h^{\alpha_0}\left(\CR\right)}\\
&\geq\left\|\left(\lambda+\left.\partial H_{\rho}\right|_{\rho=0}\right)u\right\|_{h^{\alpha_0}\left(\CR\right)}-\left\|-\frac{n-1}{R^2}\fint_{\CR}u\,d\mu_0\right\|_{h^{\alpha_0}\left(\CR\right)}\\
&\geq\frac{|\lambda-\omega|}{M}\|u\|_{h^{2+\alpha_0}_{\frac{\partial}{\partial z}}\left(\CR\right)}-\frac{n-1}{R^2}\|u\|_{C(\CR)}\\
&\geq\frac{|\lambda-\omega|}{M}\|u\|_{h^{2+\alpha_0}_{\frac{\partial}{\partial z}}\left(\CR\right)}-\frac{n-1}{R^2}\|u\|_{h^{2+\alpha_0}_{\frac{\partial}{\partial z}}\left(\CR\right)}\\
&=|\lambda|\left(\frac{\left|1-\frac{\omega}{\lambda}\right|}{M}-\frac{n-1}{R^2|\lambda|}\right)\|u\|_{h^{2+\alpha_0}_{\frac{\partial}{\partial z}}\left(\CR\right)}.
\end{align*}
By the reverse triangle inequality and by taking $Re(\lambda)\geq\omega_1:=\left(2|\omega|+\frac{2M(n-1)}{R^2}\right)$ we get
\begin{align*}
\left\|(\lambda-A_0)u\right\|_{h^{\alpha_0}\left(\CR\right)}&\geq|\lambda|\left(\frac{1}{M}-\frac{|\omega|}{M|\lambda|}-\frac{n-1}{R^2|\lambda|}\right)\|u\|_{h^{2+\alpha_0}_{\frac{\partial}{\partial z}}\left(\CR\right)}\\
&\geq\frac{|\lambda|}{2M}\|u\|_{h^{2+\alpha_0}_{\frac{\partial}{\partial z}}\left(\CR\right)},
\end{align*}
Thus, by applying Proposition 2.1.11 in \cite{Lunardi95}, $A_0$ is sectorial. This then implies by Proposition 2.4.2 in \cite{Lunardi95} that $A_{\rho}:h_{\frac{\partial}{\partial z}}^{2+\alpha_0}\left(\CR\right)\rightarrow h^{\alpha_0}\left(\CR\right)$ is sectorial for all $\rho$ in a neighbourhood of zero, $O_2\subset h^{2+\alpha_0}_{\frac{\partial}{\partial z}}\left(\CR\right)$. Since $\partial G(\rho)$ is the part of $A_{\rho}$ in $h^{\alpha}\left(\CR\right)$ the conditions of Theorem 9.2.2 in \cite{Lunardi95} are satisfied by setting $O=O_1\cap O_2$.
$\Box$

We now set
\begin{equation*}
\mathcal{S}:=\left\{\rho\in h^{2+\alpha}_{\frac{\partial}{\partial z}}\left(\CR\right):\textrm{graph(}\rho\textrm{) is a cylinder}\right\}.
\end{equation*}

\begin{lemma}
$\mathcal{M}^c$ coincides with the set $\mathcal{S}$ on a neighbourhood of zero, provided $R>\frac{d\sqrt{n-1}}{\pi}$.
\end{lemma}

\proof
Firstly, note that due to Theorem 2.4 in \cite{Simonett95}, Theorem 2.3 in \cite{Iooss92} can be applied to conclude $\mathcal{M}^c$ contains all equilibria of (\ref{GraphLin}) with $P\rho_0\in\Lambda$. Also note, the center manifold, $\mathcal{M}^c$,  is defined differently in \cite{Iooss92} as compared to \cite{Lunardi95}, however they can be seen to be equal on $\Lambda$. The rest of the proof follows from \cite{Escher98A}, where $\rho\in\mathcal{S}$ is formulated in terms of eigenfunctions:
\begin{equation*}
\rho(\bm{y})=\sum_{p=1}^{n}y_kv_{1,p,0}-R+\sqrt{\left(\sum_{p=1}^{n}y_kv_{1,p,0}\right)^2+(R+y_0)^2-\sum_{p=1}^{n}y_k^2},
\end{equation*}
with $(y_1,\ldots,y_n,0)\in\Rnone$ being the point of intersection of the axis of rotation with the $x_{n+1}$-plane and $y_0:=R'-R$ the difference between its radius and the radius of $\CR$. This map is smooth on a neighbourhood $U$ of $0\in\mathbb{R}^{n+1}$ and its derivative at zero is given by
\begin{equation*}
\partial\rho(0)\bm{y}=\sum_{p=0}^{n}y_kv_{1,p,0},\ \ \ \bm{y}\in\mathbb{R}^{n+1}.
\end{equation*}
The map taking $\bm{y}$ to the coordinates of $P\rho(\bm{y})$ with respect to the basis $v_{1,p,0}$, $0\leq p\leq n$, is then found to have derivative at zero equal to the identity and hence is a diffeomorphism from $U$ onto its image, possibly making $U$ smaller. This means that the projection of $\mathcal{S}|_{U}:=\{\rho(\bm{y}):\bm{y}\in U\}$ is an open neighbourhood of $0\in X^c$. This can be made to coincide with $\Lambda$ (after possible renaming) and since $\mathcal{S}|_{U}\subset\mathcal{M}^c$ we conclude that $\mathcal{S}$ and $\mathcal{M}^c$ coincide locally.
$\Box$

We now prove the main result.

\proof[Proof of Theorem \ref{MainR}]
As in \cite{Escher98A} we can use the same arguments in \cite{Escher98B} to prove the following result. Given $k\in\mathbb{N}$ and $\omega\in(0,-\lambda_1)$, there exists a neighbourhood $U_1(k,\omega)$ of 0 in $h^{1+\beta}_{\frac{\partial}{\partial z}}\left(\CR\right)$ with the property: Given $\rho_0\in U_1$, the solution $\rho(\cdot,\rho_0)$ of (\ref{GraphODE}) exists globally and there exist $c(k,\omega)>0$, $T(k,\omega)>0$, and a unique $z_0(\rho_0)\in\Lambda$ such that
\begin{equation*}
\|\rho(t,\rho_0)-(z_0+\gamma(z_0))\|_{C^k\left(\CR\right)}\leq ce^{-\omega t}\|(I-P)\rho_0-\gamma(P\rho_0)\|_{h^{1+\beta}\left(\CR\right)},
\end{equation*}
for $t>T$; here $z_0$ and $\gamma(z_0)$ are added in $C^k\left(\CR\right)$. This proves that $\rho(t)$ converges to an element of $\mathcal{M}^c|_{\Lambda}$, which by the above Lemma is a cylinder.
$\Box$

\begin{corollary}
Let $\Omega_0$ be a graph over a cylinder, of radius $R>\frac{d\sqrt{n-1}}{\pi}$, with height function $\rho_0$ such that the solution, $\rho(t)$, to the flow (\ref{GraphFlow}) with initial condition $\rho_0$ exists for all time and converges to zero. Then there exists a neighbourhood, $O$, of $\rho_0$ in $h_{\frac{\partial}{\partial z}}^{1+\beta}\left(\CR\right)$, $0<\beta<1$, such that for every $u_0\in O$ the solution to (\ref{GraphFlow}) with initial condition $u_0$ exists for all time and converges to a function near zero whose graph is a cylinder.
\end{corollary}

\proof
This follows by the same arguments given in \cite{Guenther02}. We set $U\subset h_{\frac{\partial}{\partial z}}^{1+\beta}\left(\CR\right)$ to be the neighbourhood of zero given in Theorem \ref{MainR}. Since $\rho(t)$ converges to zero in the $C^k$ topology, there exists a time $T$ such that $\rho(T)\in U$ and as $U$ is open there exists an open ball $B_{\epsilon}(\rho(T))\subset U$ of radius $\epsilon$ centred at $\rho(T)$. By the finite time continuity of the volume preserving mean curvature flow there exists a ball $B_{\delta}(\rho_0)$ such that if $u_0\in B_{\delta}(\rho_0)$ then the solution, $u(t)$, to (\ref{GraphFlow}) with initial condition $u_0$ exists for $t\in[0,T]$ and $u(T)\in B_{\epsilon}(\rho(T))$. Since $u(T)$ is in $U$, by Theorem \ref{MainR} the solution to (\ref{GraphFlow}) with initial condition $u(T)$ converges to a function near zero that defines a cylinder. By uniqueness of the flow we get the result.
$\Box$

\bibliographystyle{plain}

\end{document}